\documentclass[12pt]{article}
\usepackage{graphicx, mathrsfs,amssymb,amsmath,amsfonts,amsbsy,amsthm,latexsym,bm,float}
\begin{document}

\newcommand{\Dh}{\hbox{\bf D}}
\newcommand{\Xh}{\hbox{\bf X}}
\newcommand{\Fh}{\hbox{\bf F}}
\newcommand{\nablaslh}{\nabla\raise2pt\hbox{\kern-9pt\slash\kern3pt}}

\title{
XFT: An Extension of the Discrete Fractional Fourier Transform\\
\vskip.5cm}
\author{Rafael G. Campos$^{(1)}$ , J. Rico-Melgoza$^{(2)}$ and Edgar Ch\'avez$^{(1)}$\\
$^{(1)}$Facultad de Ciencias F\'{\i}sico-Matem\'aticas, \\
$^{(2)}$Facultad de Ingenier\'{\i}a El\'ectrica\\
Universidad Michoacana, \\
58060, Morelia, Mich., M\'exico.\\
\hbox{\small rcampos@umich.mx, jerico@umich.mx, elchavez@umich.mx}\\
}
\date{}
\maketitle
{
}\\
\vspace*{1truecm}
\begin{center} Abstract \end{center}
In recent years there has been a growing interest in the fractional Fourier transform driven by its large number of applications. The literature in this field follows two main routes. On the one hand, the areas where the ordinary Fourier transform has been applied are being revisited to use this intermediate time-frequency representation of signals, and on the other hand, fast algorithms for numerical computation of the fractional Fourier transform are devised. In this paper we derive a Gaussian-like quadrature of the continuous fractional Fourier transform. This quadrature is given in terms of the Hermite polynomials and their zeros. By using some asymptotic formulas, we rewrite the quadrature as a chirp-{\tt fft}-chirp transformation, yielding a fast discretization of the fractional Fourier transform and its inverse in closed form. We extend the range of the fractional Fourier transform by considering arbitrary complex values inside the unitary circle and not only at the boundary. We find that, the chirp-{\tt fft}-chirp transformation evaluated at $z=i$,  becomes a more accurate version of the {\tt fft} which can be used for non-periodic functions.

\newpage
\section{Introduction}\label{intro}
The fractional Fourier transform \cite{Nam80, Oza01}, has recently been object of renowned interest in the areas where the ordinary Fourier transform has been traditionally used. This transform can be tracked down to Wiener's paper in 1929 \cite{Wie29}. It is a particular case of the Linear Canonical Transform, which was derived in connection with canonical transformations in Quantum Mechanics \cite{Wol79, Mos71}. A standard framework to understand its properties and applications has established in \cite{Oza01}.\\
In almost any domain where the usual Fourier transform is used, there is room for techniques based on the space-frequency analysis \cite{Oza96}--\cite{Vij06}. Since this transform is a useful tool for signal processing, the direct computation of the fractional Fourier transform in digital computers has become an important issue. A survey of the different discrete implementations of this transform can be found in \cite{Tao08}. In particular some fast algorithms ${\mathcal O}(N \log N)$ have been implemented \cite{Oza96, Den97} and a comparison for the digital computation of the fractional Fourier transform has been reported in \cite{Bul04}.\\
In this paper we follow a different route by exploiting properties of the Hermite polynomials. We obtain a discrete fractional Fourier transform in terms of the standard formulation of the discrete Fourier transform. We obtain finite-dimensional vectors representing a Hermite function and its Fourier transform that converge to their exact continuous values when $N$ goes to infinity and we propose a matrix operator as a representation of the kernel of the fractional Fourier transform for functions other than the Hermite ones. This matrix, which plays the role of the kernel, yields a quadrature formula \cite{Cam92, Cam95, Cam08} and a discrete form for the fractional Fourier transform. By using some asymptotic properties of the Hermite polynomials, this discrete fractional Fourier transform can be written in terms of the standard discrete Fourier transform through a diagonal congruence transformation (a chirp-{\tt fft}-chirp transformation). In this way, an efficient algorithm of ${\mathcal O}(N\log N)$ complexity for fast computations of the fractional Fourier transform is obtained. Due to the diagonal congruence transformations, this fast fractional Fourier transform can be understood as a complex-windowed {\tt fft}. Besides, it can be evaluated in closed form in terms of exponentials for any complex value $z$ of the unit circle $\vert z\vert\le 1$ and not only at values lying on the boundary $\vert z\vert= 1$ as it is usually considered.\\
It should be said that the idea of using finite representations of the Hermite functions has widely been used to obtain discrete formulas for Fourier transforms \cite{Bul04}-\cite{Cam92}. Some of these formulas have been tested in \cite{Bul04}. The code used to test these formulas is freely available in \cite{Bul08}.\\
The rest of this paper is organized as follows. In Sec. \ref{secdos} a quadrature formula for the fractional Fourier transform is obtained. In Sec. \ref{sectres} we give the algorithm to compute a fast fractional Fourier transform and some examples are also given. Finally, a note on the convergence of this quadrature formula to the fractional Fourier transform and some remarks are given in Sec. \ref{seccuatro}.
\section{Quadrature of the fractional Fourier transform}\label{secdos}
In previous work \cite{Cam92, Cam95, Cam08} we derived a quadrature formula for the continuous Fourier transform which yields an accurate discrete Fourier transform. However, the sampling rate is far from common since the function to transform should be evaluated at the zeros of the Hermite polynomials. Now, in this section we show how to use the density properties of the Hermite polynomials to obtain both a uniform sampling rate and a quadrature of the fractional Fourier transform, yielding a discrete fractional Fourier transform.\\
Let us consider the family of Hermite polynomials $H_n(t)$, $n=0,1,\ldots$, which satisfies the recurrence equation 
\begin{equation}\label{receqg}
H_{n+1}(t)+2nH_{n-1}(t)=2tH_n(t),
\end{equation}
with $H_{-1}(t)\equiv 0$. As it is well-known \cite{Sze75}, from (\ref{receqg}) follows the Christoffel-Darboux formula \cite{Sze75}
\begin{equation}\label{chrsdar}
\sum_{n=0}^{N-1} \frac{H_n(x)H_n(y)}{2^n n!}=\begin{cases} \displaystyle\frac{H_N(x)H_{N-1}(y)-H_{N-1}(x)H_N(y)}{2^N (N-1)!(x-y)},&
x\ne y,\\ \noalign{\vskip .3cm}
\displaystyle\frac{H'_N(x)H_{N-1}(x)-H'_{N-1}(x)H_N(x)}{2^N (N-1)!},& x= y.\end{cases}
\end{equation}
Note that the recurrence equation (\ref{receqg}) can be written as the eigenvalue problem 
\begin{equation}\label{eighinf}
\begin{pmatrix} 0&1/2&0&\cdots \\1&0& 1/2&\cdots \\0& 2& 0&\cdots \\\vdots&\vdots&\vdots&\ddots\\ \end{pmatrix}
\begin{pmatrix} H_0(t)\\  H_1(t)\\  H_2(t)\\ \vdots\end{pmatrix}=t\begin{pmatrix} H_0(t)\\  H_1(t)\\  H_2(t)\\ \vdots\end{pmatrix}.
\end{equation}
Let us now consider the eigenproblem associated to the principal submatrix of dimension $N$ of (\ref{eighinf})
\[
{\mathcal H}=\begin{pmatrix}0&1/2&0&\cdots & 0& 0\\1& 0& 1/2&\cdots & 0& 0\\0& 2& 0&\cdots &  0& 0\\
\vdots&\vdots&\vdots&\ddots&\vdots&\vdots\\
0& 0& 0&\cdots & 0&1/2\\0& 0& 0&\cdots & N-1&0\end{pmatrix}.
\]
This is a general technique \cite{Gol69, Gau99} to yield gaussian quadratures for orthogonal polynomials: the recurrence equation is rewritten in matrix form to obtain orthonormal vectors of ${\mathbb R}^N$ whose entries are given (in our case) in terms of the values of $H_k(x)$, $k=0,1,\ldots,N-1$, at the zeros of $H_N(x)$. We proceed by taking a similarity transformation to symmetrize ${\mathcal H}$. 
Note that the diagonal matrix given by 
\[
S=\text{diag}\left\{1,\frac{1}{\sqrt{2}},\ldots,\frac{1}{\sqrt{(N-1)!\,2^{N-1}}}\right\},
\]
generates the symmetric matrix 
\[
H=S{\mathcal H} S^{-1}=\begin{pmatrix}0&\sqrt{\frac{1}{2}}&0&\cdots & 0& 0\\\sqrt{\frac{1}{2}}& 0& \sqrt{\frac{2}{2}}&\cdots & 0& 0\\
0& \sqrt{\frac{2}{2}}& 0&\cdots &  0& 0\\
\vdots&\vdots&\vdots&\ddots&\vdots&\vdots\\
0& 0& 0&\cdots & 0&\sqrt{\frac{N-1}{2}}\\0& 0& 0&\cdots & \sqrt{\frac{N-1}{2}}&0\end{pmatrix}.
\]
The recurrence equation (\ref{receqg}) and formula (\ref{chrsdar}) can be used to solve the eigenproblem 
\[
Hu_k=t_k u_k,\quad k=1,2,\ldots, N,
\]
which is a finite-dimensional version of (\ref{eighinf}). The eigenvalues $t_k$ are the zeros of $H_N(t)$ and the $k$th eigenvector $u_k$ is given by 
\[
c_k\left(s_1 H_0(t_k),s_2 H_1(t_k),s_3 H_2(t_k),\cdots,s_N H_{N-1}(t_k)\right)^T,
\]
where $s_1,\ldots,s_N$ are the diagonal elements of $S$ and $c_k$ is a normalization constant that can be determined from the condition $u_k^T\,u_k=1$, i.e., from
\[
c_k^2\,\sum_{n=0}^{N-1} \frac{H_n(t_k)H_n(t_k)}{2^n n!}=1.
\]
Since $H_N(t_k)=0$ and $H'_N(t_k)=2N H_{N-1}(t_k)$, the use of (\ref{chrsdar}) yields 
\[
c_k=\sqrt{\frac{2^{N-1}\,(N-1)!}{N}}\,\frac{1}{\vert H_{N-1}(t_k)\vert}=\sqrt{\frac{2^{N-1}\,(N-1)!}{N}}\,\frac{(-1)^{N+k}}{ H_{N-1}(t_k)},
\]
where we have used the fact that $\vert H_{N-1}(t_k)\vert=(-1)^{N+k} H_{N-1}(t_k)$. Thus, the components of the orthonormal vectors $u_k$, $k=1,2,\ldots, N$, are
\begin {equation}\label{ortvec}
(u_k)_n=(-1)^{N+k}\sqrt{\frac{2^{N-n}\,(N-1)!}{N\,(n-1)!}}\,\frac{H_{n-1}(t_k)}{H_{N-1}(t_k)},\quad n=1,\ldots,N.
\end{equation}
Let $U$ be the orthogonal matrix whose $k$th column is $u_k$ and $D(z)$ be the diagonal matrix $D(z)=\text{diag}\{1,z,z^2,\ldots,z^{N-1}\}$, where $z$ is a complex number. Now, let us define the matrix  
\[
{\mathcal F}_z=\sqrt{2\pi}U^{-1}D(z)U
\]
whose components are given by
\begin{equation}\label{tmat}
({\mathcal F}_z)_{jk}=\sqrt{2\pi}\,\frac{(-1)^{j+k}\,2^{N-1}\,(N-1)!}{N\,H_{N-1}(t_j)H_{N-1}(t_k)}\sum_{n=0}^{N-1}\frac{z^n}{2^n\,n!}H_n(t_j)H_n(t_k).
\end{equation}
This is the matrix representing the kernel of the fractional Fourier transform in a $N$-dimensional vector space, as we show next.
\subsubsection*{Asymptotic formulae}
Let us look for the asymptotic form of the components (\ref{tmat}) of ${\mathcal F}_z$. First note that the asymptotic expression for $H_N(t)$ in the oscillatory region is (Eq. (8.22.8) of \cite{Sze75})
\begin{equation}\label{asymhnt}
H_N(t)=\frac{\Gamma(N+1)}{\Gamma(N/2+1)}e^{t^2/2}\left(\cos(\sqrt{2N+1}\,\,t-N\pi/2)+{\cal O}(N^{-1/2})\right).
\end{equation}
This gives an approximate form for the zeros of $H_N(t)$, i.e.,
\begin{equation}\label{asymcer}
t_k= \left(\frac{2 k-N-1}{\sqrt{2N}}\right)\frac{\pi}{2},
\end{equation}
$k=1,2,\ldots,N$. Thus, the use of (\ref{asymhnt}) and (\ref{asymcer}) yields
\[
H_{N-1}(t_k)\simeq (-1)^{N+k}\,\frac{\Gamma(N)}{\Gamma(\frac{N+1}{2})}\,e^{t_k^2/2},\quad N\to\infty.
\]
Therefore, for $N$ large enough, (\ref{tmat}) can be written as
\[
({\mathcal F}_z)_{jk}\simeq \sqrt{2\pi}\,\frac{2^{N-1}\,[\Gamma(\frac{N+1}{2})]^2}{\Gamma(N+1)}\,\,e^{-(t_j^2+t_k^2)/2}\,
\sum_{n=0}^{\infty}\frac{z^n}{2^n\,n!}H_n(t_j)H_n(t_k),
\]
and  finally, Stirling's formula and Mehler's formula \cite{Erd53} produce the result
\begin{equation}\label{fasy}
({\mathcal F}_z)_{jk}\simeq \sqrt{\frac{2}{1-z^2}}\,\exp\left(-\frac{(1+z^2)(t_j^2+t_k^2)-4t_jt_kz}{2(1-z^2)}\right)\Delta t_k,
\end{equation}
where $\Delta t_k$ is the difference between two consecutive asymptotic Hermite zeros, i.e.,
\begin{equation}\label{deltk}
\Delta t_k=t_{k+1}-t_k=\frac{\pi}{\sqrt{2N}}.
\end{equation}
Let us consider now a complex-valued function $g(t)$ defined for $t\in{\mathbb R}$ and let us form the vector
\[
g=(g(t_1),g(t_2),\ldots,g(t_N))^T.
\] 
Therefore, the multiplication of the matrix ${\mathcal F}_z$ by the vector $g$ gives the vector $G$ with entries
\begin{equation}\label{primcuad}
G_j=\sum_{k=1}^N({\mathcal F}_z)_{jk}g(t_k)\simeq \sqrt{\frac{2}{1-z^2}}\,\sum_{k=1}^N\exp\left(-\frac{(1+z^2)(t_j^2+t_k^2)-4t_jt_kz}{2(1-z^2)}\right)g(t_k)\Delta t_k,
\end{equation}
for $j=1,2\ldots,N.$ Note that this equation is a Riemann sum for the integral
\begin{equation}\label{eqifft}
{\mathscr F}_z[g(t'),t]=\sqrt{\frac{2}{1-z^2}}\,\int_{-\infty}^\infty \exp\left(-\frac{(1+z^2)(t^2+\mathop{t'}^2)-4t\mathop{t'}z}{2(1-z^2)}\right)g(t')dt', \quad \vert z\vert<1,
\end{equation}
that is,
\begin{equation}\label{cuadtrfor}
{\mathscr F}_z[g(t'),t_j]\simeq\sum_{k=1}^N({\mathcal F}_z)_{jk}g(t_k),\quad N\to\infty.
\end{equation}
Note that ${\mathscr F}_z[g(t'),t]$ is the continuous fractional Fourier transform \cite{Nam80} of $g(t')$ up to a constant. Thus, the matrix ${\mathcal F}_z$ is a discrete fractional Fourier transform.  It should be noted that it is a discretization of the fractional Fourier transform for any complex value $z$ of the unitary circle $\vert z\vert\le 1$ and not only for $z$ lying on its boundary, as it is usually considered and that the argument $\varphi$ of $z=r\exp(i\varphi)$ is real and not complex as it is used in some applications \cite{Pel06}. \\
In the case $z=\pm i$, (\ref{fasy}) becomes
\[
F_{jk}\equiv ({\mathcal F}_{\pm i})_{jk}\simeq e^{\pm it_jt_k}\Delta t_k,
\]
thus, we have a discrete Fourier transform. It has been previously obtained in \cite{Cam92} where some numerical examples are given,  and it has been applied to the analysis of brain signals \cite{Sol09}. 
\section{A fast discrete fractional Fourier transform}\label{sectres}
The {\tt fft} can be used to obtain a fast algorithm for (\ref{cuadtrfor}) as follows. Since the matrix ${\mathcal F}_z$ represents a quadrature for the fractional Fourier transform expected to converge when a large number of nodes are used, we consider the asymptotic form (\ref{fasy}) which can be written in more detail as
\begin{equation}\label{fasydos}
({\mathcal F}_z)_{jk}=\sqrt{\frac{2}{1-z^2}}\,\exp(-\mu t_j^2)\exp(\nu t_jt_k)\exp(-\mu t_k^2)\Delta t_k,
\end{equation}
where we have used the definitions
\begin{equation}\label{munufffrt}
\mu=\frac{1+z^2}{2(1-z^2)}, \quad \nu=\frac{2z}{1-z^2}.
\end{equation}
Note that we have replaced the approximately equal sign "$\simeq$" by the equal sign "$=$" in (\ref{fasy}) redefining ${\mathcal F}_z$ in (\ref{fasydos}). To distinguish the fast implementation of the present discrete fractional transform from other important contributions, we denote it by {\tt xft} (a extended Fourier transform). To show the main differences of this {\tt xft} with the usual {\tt fft}, we consider first the case of the standard Fourier transform.
\subsection{The {\tt xft} as an improvement of the {\tt fft}}
As noted above, the case $z=i$ in (\ref{fasydos}) corresponds to a discrete Fourier transform 
\begin{equation}\label{fftuno}
F_{jk}=\frac{\pi}{\sqrt{2N}}\exp\left[i\frac{\pi^2}{2N} \left(j-\frac{N-1}{2}\right) \left(k-\frac{N-1}{2}\right)\right]
\end{equation}
where now $j,k=0,1,2,\ldots,N-1$, and we have used (\ref{asymcer}) and (\ref{deltk}). Since $\sum_{k=1}^NF_{jk}g(t_k)$ is a quadrature and therefore, an approximation [cf. Eq. \ref{cuadtrfor}] of
\[
G(\omega_j)=\int_{-\infty}^\infty e^{i \omega_j t} g(t)dt,
\]
we can use the basic property of the Fourier transform
\begin{equation}\label{furesc}
G(a \omega_j)=\int_{-\infty}^\infty e^{i a \omega_j t} g(t)dt
\end{equation}
to rewrite (\ref{fftuno}) in the more convenient form
\begin{equation}\label{fftdos}
(F_a)_{jk}=\frac{\pi}{\sqrt{2N}}\exp\left[ia\frac{\pi^2}{2N} \left(j-\frac{N-1}{2}\right) \left(k-\frac{N-1}{2}\right)\right],
\end{equation}
in the understanding that this matrix yields a scaled discrete Fourier transform. The scaling parameter $a$ can be used to connect the {\tt xft} with the standard form of the discrete Fourier transform. If we choose $a=4/\pi$, Eq. (\ref{fftdos}) can be written as
\[
(F_{4/\pi})_{jk}=\frac{\pi e^{i\frac{\pi}{2}\frac{(N-1)^2}{N}}}{\sqrt{2N}}\left[e^{-i\pi\frac{N-1}{N} j}\right]\left[e^{i\frac{2\pi}{N}jk}\right]
\left[e^{-i\pi\frac{N-1}{N} k}\right],
\]
where $ j,k=0,1,2,\ldots,N-1$.  In matrix form, and defining $\tilde{F}\equiv F_{4/\pi}$ to simplify the notation, the following discrete Fourier transform $\tilde{F}$ can be given in terms of the usual discrete Fourier transform $D_F$ as
\begin{equation}\label{ffttres}
\tilde{F}=\frac{\pi e^{i\frac{\pi}{2}\frac{(N-1)^2}{N}}}{\sqrt{2N}}  SD_FS,
\end{equation}
where $S$ is the diagonal matrix whose nonzero elements are $\exp(-i\pi (N-1) j/N)$, $j=0,1,\ldots,N-1$ and $D_F(\cdot)$ can be computed through a standard {\tt fft} algorithm. Therefore, the computational cost of $\tilde{F}(\cdot)$ is the same of the algorithm used to compute the {\tt fft} of the vector $S(\cdot)$. However, the output of the present method is more precise since it comes from a convergent quadrature formula. This is shown below.\\
The inverse of $\tilde{F}$ is easily obtained:
\[
(\tilde{F}^{-1})_{jk}=\frac{\sqrt{2/N}}{\pi}\exp\left[-i\frac{2\pi}{N} \left(j-\frac{N-1}{2}\right) \left(k-\frac{N-1}{2}\right)\right],
\]
where $ j,k=0,1,2,\ldots,N-1$.\\
In the applications, we have to remind that $\tilde{F}$ gives a scaled transform. The following simple algorithm incorporate these ideas.\\

\fbox{
\begin{minipage}{14cm}
\begin{center} \vskip 0.4cm
Algorithm 1 \end{center}
\hrule width 14cm \vskip.3truecm
To compute an approximation $G=(G_1,G_2,\ldots,G_N)^T$ of the Fourier transform (\ref{furesc}) of the vector
\[
g=(g_1,g_2,\ldots,g_N)^T.
\]
\begin{enumerate}
\item For given $N$ set $t_k=\pi (2 k-N-1)/[2(\sqrt{2N})]$,\,\, $k=1,2,\ldots,N.$
\item Let $D_F$ be the discrete Fourier transform, i.e., $(D_F)_{jk}=e^{i\frac{2\pi}{N}jk}$. For $j,k=0,1,2,\ldots,N-1$, compute the diagonal matrix $S$ according to $S_{jk}=e^{-i\pi\frac{N-1}{N} j}\delta_{jk}$.
\item Obtain the approximation $G_j$ to $G(\frac{4}{\pi} t_j)$ by computing the matrix-vector product
\begin{equation}\label{algfxft1}
G=\frac{\pi e^{i\frac{\pi}{2}\frac{(N-1)^2}{N}}}{\sqrt{2N}} SD_F(Sg),
\end{equation}
with a standard {\tt fft} algorithm.
\end{enumerate}
\end{minipage}}

\noindent
Either if the input vector $g$ is given by the values of a function $g(t)$ at $t_k$, or not, we can represent  the approximation to the Fourier transform $G(\omega)$, as given in a plot of points $(\frac{4}{\pi} t_j,G_j)$.\\
In the examples given below we plot the real and imaginary parts of the vector $G$ and the exact transform $G(\frac{4}{\pi} \omega)$. We give first two examples of non-periodic/singular functions \cite{Erd53b}. We compare the performance of the scaled {\tt xft} evaluated at $z=i$ with that of the standard {\tt fft}.
\subsubsection*{Example 1}
Consider the pair of Fourier transforms
\[
g(t)=\cos(t^2),\qquad G(\omega)=\sqrt{\pi}\cos(\frac{\omega^2-\pi}{4}).
\]
Figures 1, 2 and 3 show the numerical convergence attained by the {\tt xft} for N=512 and 1024 respectively and it is compared with the {\tt fft} output. The plots on the left-hand side of Figs. 1 and 2, show the output of {\tt xft} (dashed lines) plotted against the exact values of $G(\frac{4}{\pi} \omega_j)$  (solid lines) with
\[
\omega_j=\pi (j-N/2-1/2)/\sqrt{2N},\quad j=1,2,\ldots,N.
\]
The plots on the right-hand side of Figs. 1 and 2, show in dashed lines, $G^{\text{fft}}$, the output of {\tt fft}, plotted against the exact values $G(\omega_j)$ (solid lines). In Fig. 3A we show the error $\vert G(\frac{4}{\pi} \omega_j)-G_j\vert$. In Fig. 3B, the error $\vert G(\omega_j)-G^{\text{fft}}_j\vert$ is shown.
\begin{figure}[H]\label{Figun}
\centering
\includegraphics[scale=0.55]{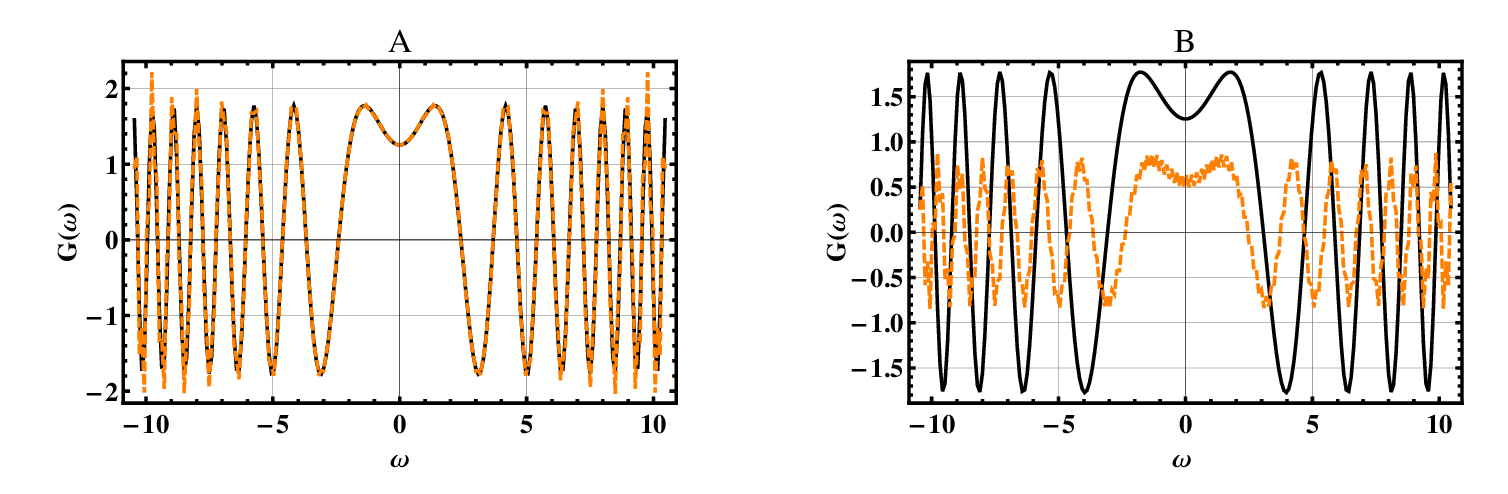}
\caption{{\bf A:} Part of the plots corresponding to $(\omega_j,G(\frac{4}{\pi}\omega_j))$ (solid line) and $(\omega_j,G_j)$ (dashed line) for $N=512$. {\bf B:} Part of the plots corresponding to $(\omega_j,G(\omega_j))$ (solid line) and $(\omega_j,G^{\text{fft}}_j)$ (dashed line) computed for the same function and number of points.}
\end{figure}

\begin{figure}[H]\label{Figdos}
\centering
\includegraphics[scale=0.55]{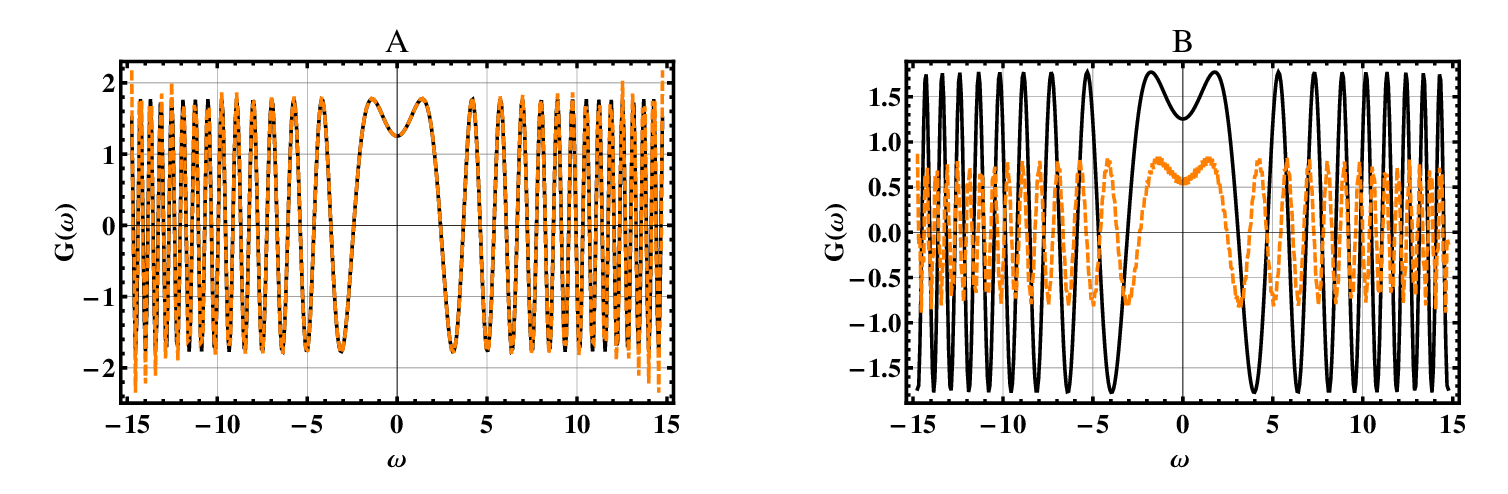}
\caption{{\bf A:} Part of the plots corresponding to $(\omega_j,G(\frac{4}{\pi}\omega_j))$ (solid line) and $(\omega_j,G_j)$ (dashed line) for $N=1024$. {\bf B:} Part of the plots corresponding to $(\omega_j,G(\omega_j))$ (solid line) and $(\omega_j,G^{\text{fft}}_j)$ (dashed line) computed for the same function and number of points.}
\end{figure}

\begin{figure}[H]\label{Figtres}
\centering
\includegraphics[scale=0.55]{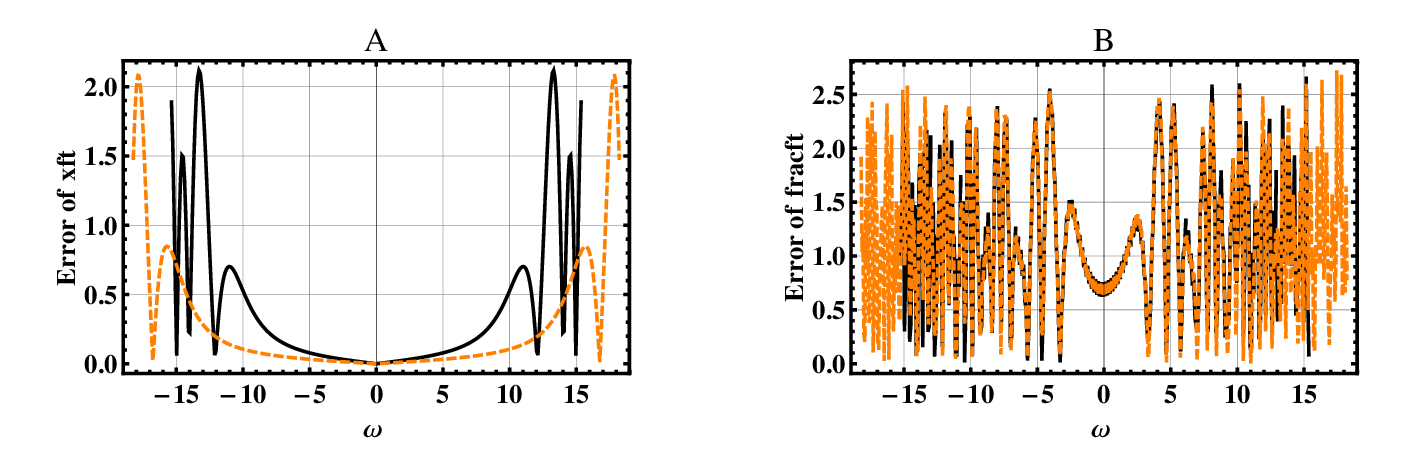}
\caption{{\bf A:} Error of the output produced by {\tt xft}, $\vert G(\frac{4}{\pi} \omega_j)-G_j\vert$, for $N=512$ (solid line) and $N=1024$ (dashed line). 
{\bf B:} Error of the output produced by {\tt fft}, $\vert G(\omega_j)-G^{\text{fft}}_j\vert$, for $N=512$ (solid line) and $N=1024$ (dashed line).}
\end{figure}

The norm $\max_{j=1}^N \vert G^{\text{exact}}_j-G^{\text{approx}}_j\vert$, evaluated for $N=512$ and $N=1024$, is 2.11 and 2.08 for {\tt xft}, and 2.68 and 2.65 for {\tt fft}, respectively.\\
It can be seen from Fig. 3A that the output of {\tt xft} converges to the exact transform. The error yielded by the {\tt xft} approaches zero at more points around the origin as $N$ become large, whereas the error produced by the {\tt fft} maintains almost the same magnitude on the same interval.
\subsubsection*{Example 2}
As an example of a singular function, let us consider the pair of Fourier transforms (the integral is a Cauchy Principal Value)
\[
g(t)=\frac{e^{-t/2}}{2-e^{-t}},\qquad G(\omega)=\frac{\pi}{2^{1/2+i \omega}}\cot(\frac{\pi}{2}-i\pi \omega).
\]
Figures 4 and 5 shows the performance of the {\tt xft} for N=512 compared with that of the {\tt fft}. We compare the exact transform $(\omega_j,G(\frac{4}{\pi}\omega_j))$ (solid line) against the output of {\tt xft} $(\omega_j,G_j)$ (dashed line) in Fig. 4 and the exact transform $(\omega_j,G(\omega_j))$ (solid line) against the output of {\tt fft} $(\omega_j,G^{\text{fft}}_j)$ (dashed line) in Fig. 5. 
\begin{figure}[H]\label{Figcuatro}
\centering
\includegraphics[scale=0.55]{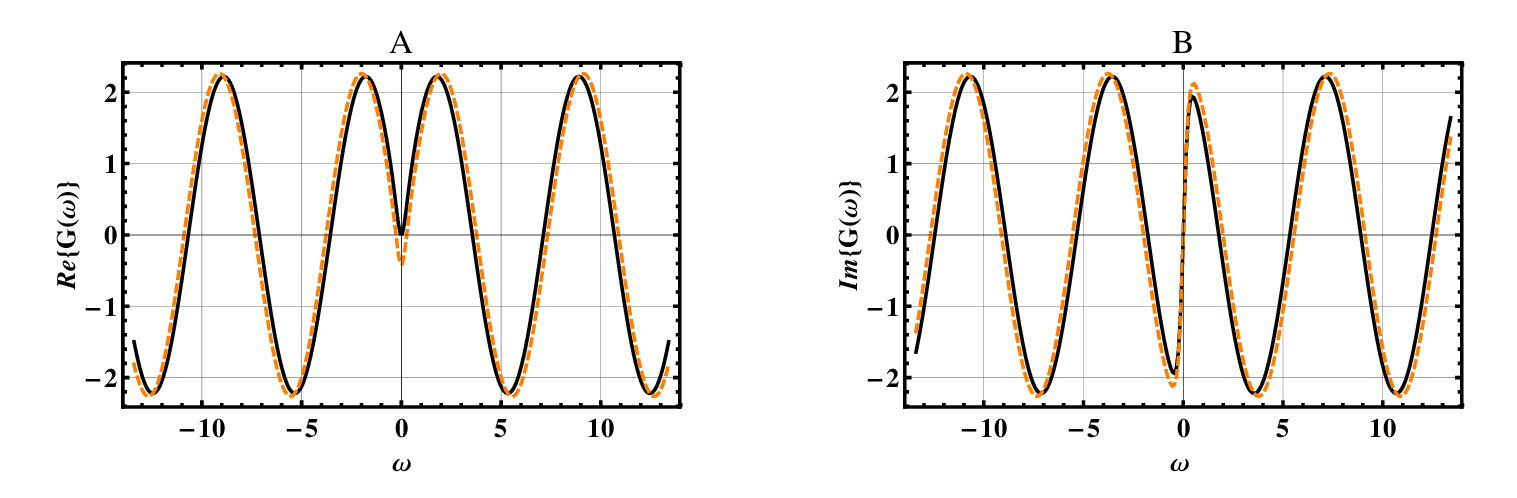}
\caption{Real part (plot {\bf A}) and imaginary part (plot {\bf B}) of the continuous Fourier transform (solid line) compared with the output of {\tt xft} (dashed line) computed with $N=512$. The function $g(t)$ is that given in Example 2. The max-norm of the error is 0.4262 for both the real and imaginary part.}
\end{figure}

\begin{figure}[H]\label{Figcinco}
\centering
\includegraphics[scale=0.55]{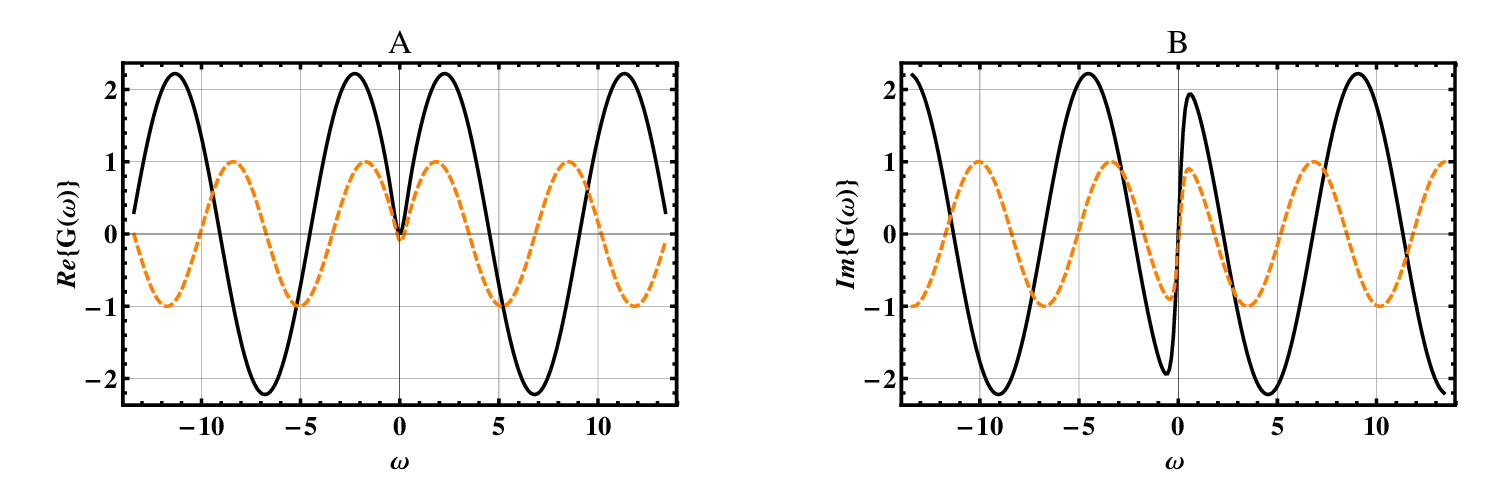}
\caption{Real part (plot {\bf A}) and imaginary part (plot {\bf B}) of the continuous Fourier transform (solid line) compared with the output of the {\tt fft} (dashed line) with $N=512$. The function $g(t)$ is that given in Example 2. The max-norm of the error is 3.18 and 3.21 for the real and imaginary part, respectively.}
\end{figure}
\subsubsection*{Example 3}
This example is concerned with the performance of the {\tt xft} acting on sinusoidal waveforms of integer-harmonic frequency.  Our starting point is again Eq. (\ref{fftdos}). By changing the order of the indexes, (\ref{fftdos}) can be written as 
\[
(\tilde{F})_{jk}=\frac{\pi}{\sqrt{2N}}e^{i\frac{2\pi}{N} jk},
\]
where
\[
j,k =\begin{cases} 0,\pm 1,\pm 2,\cdots,\pm (N-1)/2,& \text{N odd},\\ \pm 1/2,\pm 3/2,\cdots,\pm (N-1)/2,& \text{N even}.\end{cases}
\]
Let $N$ be odd, and assume that the signal to be transformed is an integer-harmonic signal. For instance, 
\begin{equation}\label{ftcosta}
g_k=\cos(k\tau_m),\quad \tau_m=\frac{2\pi m}{N},
\end{equation}
where $k=0,\pm 1,\pm 2,\ldots,\pm (N-1)/2.$ Then, the formulae \cite{Obe73} 
\[
\sum_{k=1}^N \cos(kx)=\frac{\sin(N x/2)\cos((N+1) x/2)}{\sin(x/2)},\quad \sum_{k=1}^N \sin(kx)=\frac{\sin(N x/2) \sin((N+1) x/2)}{\sin(x/2)},
\]
can be used to prove that the {\tt xft} transform of (\ref{ftcosta}) is
\begin{equation}\label{xftinthar}
G_j=\sum_{k=-(N-1)/2}^{(N-1)/2} (\tilde{F})_{jk}\cos(k\tau_m)=\pi\sqrt{\frac{N}{2}} (\delta_{j,-m}+\delta_{jm}), 
\end{equation}
where $j=0,\pm 1,\pm 2,\cdots,\pm (N-1)/2$, corresponding to two pulses centered at the integers $m$ and $-m$ respectively.\\
Now assume that $N$ is even and take a half-integer harmonic signal, i.e., one of the same form as (\ref{ftcosta}) but with $k=\pm 1/2,\pm 3/2,\ldots,\pm (N-1)/2$. A similar calculation as above shows that $G_j$ has the same form of (\ref{xftinthar}) but the sum runs over half-integers and the pulses are now centered at the half-integers $m$ and $-m$ respectively.\\
In the two above cases the dispersion is zero, but this does not happen if we take $N$ odd and consider a half-integer harmonic signal, or $N$ even and consider an integer harmonic signal. Finally note that these processes can be reversed: the {\tt xft} of the corresponding pulses yields integer or half-integer harmonic signals.\\
We give another numerical case with a cosine waveform with an arbitrary frequency. Let us consider the example $g(t)=\cos(5.156\,t)$. In Fig. 6 we show the {\tt xft} output for ($a$) $N=1024$ and ($b$) $N=2048$. As it can be seen, the leakage practically vanishes for the latter case. As a measure of the leakage we give the mean value
\[
\mu=\frac{1}{N}\sum_{k=1}^N (\vert G_k\vert - \delta_{km}\vert G_m\vert - \delta_{km'}\vert G_{m'}\vert),
\]
where the indexes $m$ and $m'$ correspond to the frequencies where $\vert G\vert$ attains their maxima. Thus, we find that $\mu =0.14105$ for $N=1024$ and $\mu =0.00276$ for $N=1024$.
\begin{figure}[H]\label{Figseis}
\centering
\includegraphics[scale=0.55]{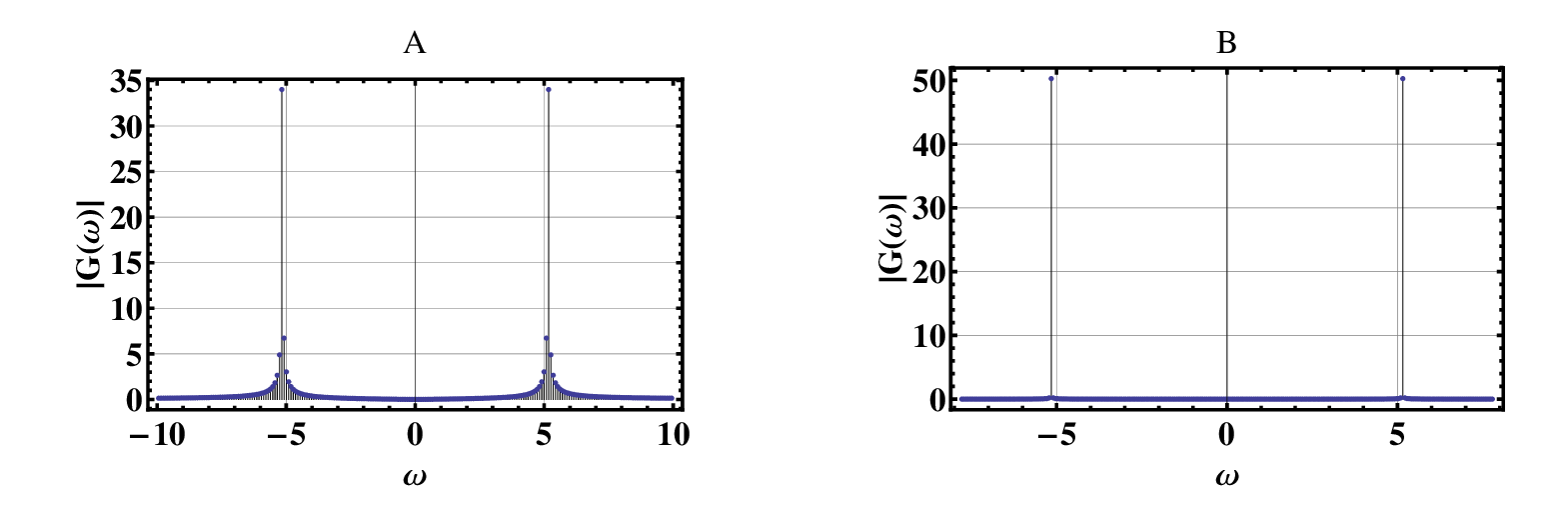}
\caption{Fourier transform of $\cos(5.156\,t)$. The maximum occurs at $\omega=5.17072$ for $N=1024$ in plot {\bf A} and at $\omega=5.15625$ for $N=2048$ in plot {\bf B}.}
\end{figure}
\subsection{The {\tt xft} as a fast discrete fractional Fourier transform}
To obtain a discrete and fast implementation of the continuous fractional Fourier transform, we follow the same procedure as above. The definition of the fractional transform given by (\ref{eqifft}) is $\sqrt{2 \pi}$ times the one given in \cite{Nam80}. This has to be taken into account for comparing purposes. Our starting point is the definition of the fractional transform given in (\ref{eqifft})
\[
G_z(t)\equiv{\mathscr F}_z[g(t'),t]=\sqrt{\frac{2}{1-z^2}}\,\, e^{-\mu t^2}\int_{-\infty}^\infty e^{ \nu t\mathop{t'}} e^{-\mu\mathop{t'}^2}g(t')dt',
\]
where we have used the definitions (\ref{munufffrt}). The simple scaling
\begin{equation}\label{aftftg}
G_z(at)=\sqrt{\frac{2}{1-z^2}}\,\, e^{-\mu a^2t^2}\int_{-\infty}^\infty e^{a \nu t\mathop{t'}} e^{-\mu\mathop{t'}^2}g(t')dt',
\end{equation}
and the discretization of the kernel
\[
({\mathcal F}^a_z)_{jk}=\sqrt{\frac{2}{1-z^2}}\,\exp(-\mu a^2t_j^2)\exp(a\nu t_jt_k)\exp(-\mu t_k^2)\Delta t_k,
\]
allows us to implement a fast algorithm to compute the {\tt xft} in terms of the {\tt fft}. This can be done by choosing $a=i2(1-z^2)/(\pi z)$, because then 
\begin{equation}\label{eqxftfrac}
({\mathcal F}^a_z)_{jk}=\sqrt{\frac{2}{1-z^2}}\,\exp(-\mu a^2t_j^2) (\tilde{F})_{jk}\exp(-\mu t_k^2),\quad a=\frac{2i}{\pi z}(1-z^2),
\end{equation}
where $\tilde{F}$, given by (\ref{ffttres}), is the fast {\tt xft} for the case of the standard Fourier transform. In order to simplify the notation, let $\tilde{\mathcal F}_z$ be ${\mathcal F}^{2i(1-z^2)/(\pi z)}_z$. Since $\sum_{k=1}^N (\tilde{\mathcal F}_z)_{jk}f(t_k)$ is a quadrature of (\ref{aftftg}) with $a=i2(1-z^2)/(\pi z)$ [cf. Eq. \ref{cuadtrfor}], $\tilde{\mathcal F}_z$ is a fast discrete fractional Fourier transform which gives an approximation for the scaled function $G_z(at)$ at the nodes $t_k$.\\
We remind the reader that the parameters $\mu$ and $\nu$ depend on the complex number $z$ [cf. Eq. (\ref{munufffrt})], and that (\ref{eqxftfrac}) gives an approximation to the fractional Fourier transform, as defined by (\ref{eqifft}), for any complex number $z$, $\vert z\vert\le 1$ and not only for $z$ lying on its boundary, as it is usually considered.\\
Note that the inverse transform can be obtained in terms of $\tilde{F}^{-1}$. These ideas are incorporated in the following

\fbox{\begin{minipage}{15cm}
\begin{center} \vskip 0.4cm
Algorithm 2 \end{center}
\hrule width 15cm \vskip.3truecm
To compute an approximation $G_z=(G_{z1},G_{z2},\ldots,G_{zN})^T$ of the fractional Fourier transform ${\mathscr F}_z[g(t'),i2(1-z^2)/(\pi z)t_j]$ of the vector 
\[
g=(g_1,g_2,\ldots,g_N)^T,
\]
at the complex value $z$, $\vert z\vert\le1$.
\begin{enumerate}
\item For given $N$ and $z$, set
\begin{enumerate}
\item $\mu=(1+z^2)/[2(1-z^2)]$,\,\, $a=2i(1-z^2)/(\pi z)$.
\item $t_k=\pi (2 k-N-1)/[2(\sqrt{2N})]$,\,\, $k=1,2,\ldots,N$
\end{enumerate}

\item Let $D_F$ be the discrete Fourier transform, i.e., $(D_F)_{jk}=e^{i\frac{2\pi}{N}jk}$. For $j,k=0,1,2,\ldots,N-1$, compute the diagonal matrices $S_1$ and $S_2$ according to 
\[(S_1)_{jk}=e^{-\mu a^2 t_j^2-i\pi\frac{N-1}{N} j}\delta_{jk},\quad
(S_2)_{jk}=e^{-\mu t_j^2-i\pi\frac{N-1}{N} j}\delta_{jk}.\]
\item Obtain the approximation $G_{zj}$ to $G_z(a t_j)$ by computing the matrix-vector product
\begin{equation}\label{algfxft2}
G_z=\sqrt{\frac{2}{1-z^2}}\,\frac{\pi e^{i\frac{\pi}{2}\frac{(N-1)^2}{N}}}{\sqrt{2N}} S_1D_F(S_2g),
\end{equation}
with a standard {\tt fft} algorithm.
\end{enumerate}
\end{minipage}}

\noindent 
To compare the performance of the {\tt xft} with other implementations we have used the software that is available in \cite{Bul08}, in particular the  {\tt fracft} routine.\\
We test this fast transform on two well-known examples for which $\vert z \vert=1$ (see Ref. \cite{Nam80}). In examples 4 and 5, $a=2i(1-z^2)/(\pi z)$ and $\omega_k=t_k$. In these examples, we compare the exact transform $(\omega_j,G(\frac{4}{\pi}\omega_j))$ (solid line) against the output of {\tt xft} $(\omega_j,G_j)$ (dashed line) in the left-hand side plots, and the exact transform $(\omega_j,G(\omega_j))$ (solid line) against the output of {\tt fracft} $(\omega_j,G^{\text{fracft}}_j)$ (dashed line) in the right-hand side plots. 

\subsubsection*{Example 4}
Consider the pair of fractional Fourier transforms
\[
g(t)=\exp(-t^2/2+\beta t),\qquad G_z(t)=\exp(-t^2/2-(i/2) \beta^2 e^{i\varphi}\sin\varphi+\beta t e^{i\varphi}),\,\,\, z=e^{i\varphi}.
\]
Fig. 7A shows a part of the plot of the real part of the {\tt xft} compared with the real part of $G_z(a\omega_j)$. The max-norm of the error is of order $10^{-12}$ for the real (and imaginary) part yielded by {\tt xft}. In Fig. 7B we show the output of {\tt fracft} for the real part of $G_z(a\omega_j)$. 
\begin{figure}[H]\label{Figsiete}
\centering
\includegraphics[scale=0.55]{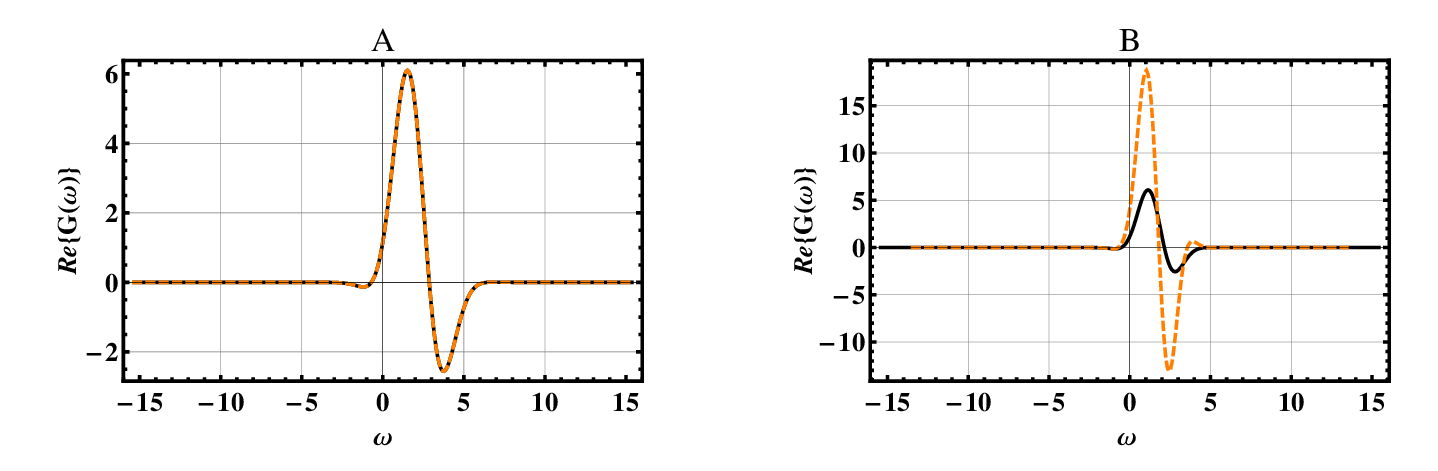}
\caption{Real part of the fractional Fourier transform (solid lines) compared with the output of {\tt xft} in {\bf A} and with the output of {\tt fracft} in {\bf B} (dashed lines), computed with $N=512$. The function $g(t)$ is that given in Example 4 for $\beta=2$ and $\varphi=\pi/5$.}
\end{figure}
\subsubsection*{Example 5}
Consider the pair of fractional Fourier transforms
\[
g(t)=1,\qquad G_z(t)=\exp[i (t^2 \tan\varphi- \varphi )/2]/\sqrt{\cos\varphi},\,\,\, z=e^{i\varphi}.
\]
Fig. 8 shows the performance of the {\tt xft} and {\tt fracft} compared with the real part of the exact transform $G_z(a\omega_j)$ in Figs. 8A and 8B respectively. 
\begin{figure}[H]\label{Figocho}
\centering
\includegraphics[scale=0.55]{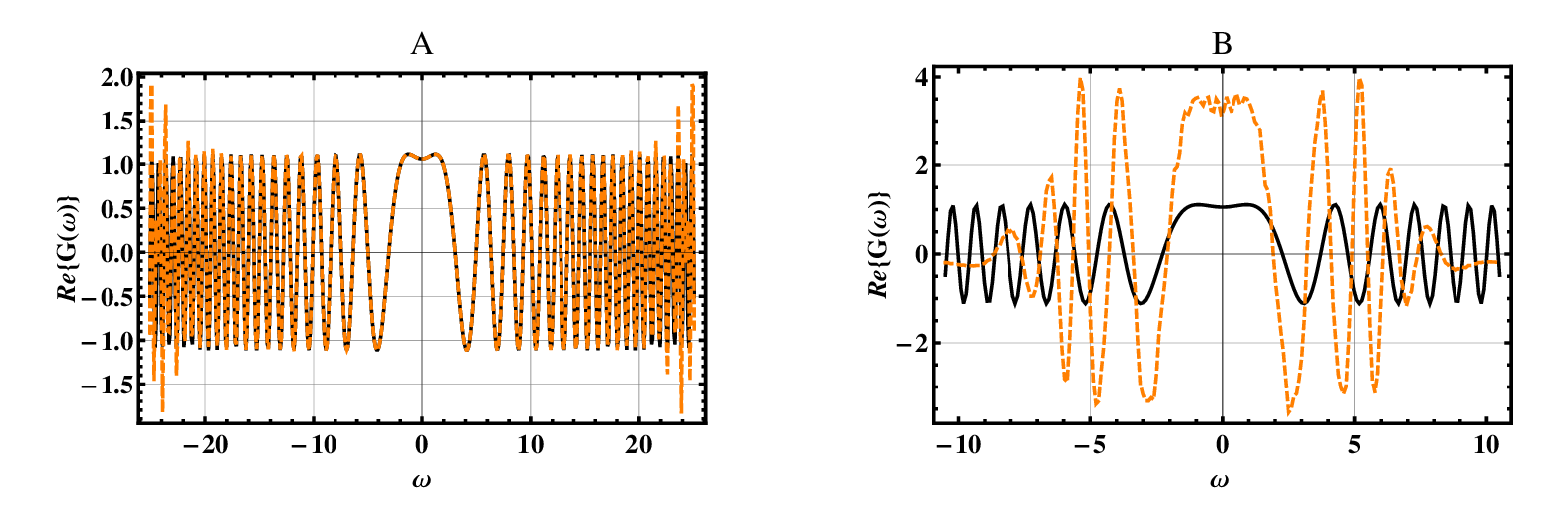}
\caption{Real part of the fractional Fourier transform (solid lines) compared with the output of {\tt xft} in {\bf A} and with the output of {\tt fracft} in {\bf B} (dashed lines), computed with $N=512$. The function $g(t)$ is that given in Example 5 for $\varphi=\pi/5$.}
\end{figure}
All the examples given in this section show that the {\tt xft} can be used to produce an accurate approximation to the fractional Fourier transform. 
\section{Conclusions and final remarks}\label{seccuatro}
Some points concerning the fast algorithm to compute the fractional Fourier transform presented in this paper have to be noticed. First, note that the {\tt xft} and their inverse have simple forms given by explicit matrix elements that can be computed easily for complex values of the parameter $z=r\exp(i\varphi)$ lying inside the unitary circle, not only on the boundary. The {\tt xft} can be interpreted as a complex-windowed {\tt fft}, as fast as the latter. The {\tt xft} performs well with singular functions or non-periodic functions and it is able to detect half-integer harmonic frequencies with no leakage.\\
Worth of note is the fact that the xft can compute a fractional Fourier transform with the complex parameter inside the unitary circle, and that this case has not been studied or used in an application. Further research is necessary to give applications of this extension. 



\begin{thebibliography}{99}
\bibitem{Nam80} Namias, V., {\sl The Fractional Order Fourier Transform and its Application to Quantum Mechanics}, J. Inst. Maths. Applics. {\bf 25} (1980) 241-265.
\bibitem{Oza01} Ozaktas, H.M; Zalevsky, Z; Kutay, M.A., {\sl The Fractional Fourier Transform with Applications in Optics and Signal Processing}, John Wiley and Sons, Chichester, UK, 2001.
\bibitem{Wie29} Wiener, N., {\sl Hermitian polynomials and Fourier analysis}, J. Math. and Phys. MIT {\bf 8} (1929) 70-73.
\bibitem{Wol79} Wolf, K.B., {\sl Integral Transforms in Science and Engineering}, Ch. 9-10, Plenum Press, New York, 1979.
\bibitem{Mos71} Moshinsky M.; Quesne, C., {\sl Linear Canonical Transformations and their unitary representations}, J. Math. Phys., {\bf 12} (1971), 1772-1783.
\bibitem{Oza96} Ozaktas, H.M.; Ankan, O,; Kutay, M.A.;  Bozda\u{g}i, G., {\sl Digital Computation of the Fractional Fourier Transform}, IEEE Trans. Signal Processing, {\bf 44} (1996) 2141-2150.
\bibitem{Sax05} Saxena R.; and Singh, K., {\sl Fractional Fourier transform: A novel tool for signal processing}, J. Indian Inst. Sci. {\bf 85}  (2005) 11-26.
\bibitem{Oon08} Oonincx, P.J.,{\sl Joint time-frequency offset detection using the fractional Fourier transform}, Signal Processing {\bf 88} (2008) 2936-2942.
\bibitem{Sha07} Sharma, K.K.; Joshi, S.D., {\sl Time delay estimation using fractional Fourier transform}, Signal Processing {\bf 87} (2007) 853-865.
\bibitem{Onu07} Onural, L.; Ozaktas, H.M., {\sl Signal processing issues in diffraction and holographic 3DTV}, Signal Processing: Image Communication {\bf 22} (2007) 169-177.
\bibitem{Vij06} Vijaya, C.; Bhat, J.S., {\sl Signal compression using discrete fractional Fourier transform and set partitioning in hierarchical tree}, Signal Processing {\bf 86} (2006) 1976-1983.
\bibitem{Tao08} Tao, R.; Zhang, F.; Wang, Y., {\sl Research progress on discretization of fractional Fourier transform}, Sci. China Ser. F {\bf 51} (2008) 859-880.
\bibitem{Den97} Deng, X.; Li, Y.; Fan, D.; Qiu, Y., {\sl A fast algorithm for fractional Fourier transform}, Opt. Commun. {\bf 138} (1997) 270-274.
\bibitem{Bul04} Bultheel, A.; Sulbaran, H.M., {\sl Computation of the fractional Fourier transform}, Appl. Comput. Harmon. Anal. {\bf 16} (2004) 182-202. 
\bibitem{Can98} Candan, \c{C}., {\sl The discrete Fractional Fourier Transform}. MS Thesis, Bilkent University, Ankara, 1998.
\bibitem{Can00} Candan, \c{C},; Kutay, M.A.; Ozaktas, H.M., {\sl The discrete fractional Fourier transform}, IEEE Trans. Sig. Proc. {\bf 48} (2000) 1329-1337
\bibitem{Can07} Candan,  \c{C}., {\sl On higher order approximations for Hermite-Gaussian functions and discrete fractional Fourier transforms}, IEEE Signal Process. Lett.  {\bf 14} (2007) 699-702.
\bibitem{Pei99} Pei,S.C.; Yeh, M.H.; Tseng, C.C., {\sl Discrete fractional Fourier transform based on orthogonal projections}, IEEE Trans. Sig. Proc. {\bf 47} (1999) 1335-1348
\bibitem{Pei06} Pei,S.C.; Hsue, W.L.; Ding, J.J., {\sl Discrete fractional Fourier transform based on new
nearly tridiagonal commuting matrices}, IEEE Trans. Sig. Proc. {\bf 54} (2006) 3815-3828
\bibitem{Cam92} Campos, R.G.; Ju\'arez, L.Z., {\sl A discretization of the Continuous Fourier Transform}, Il Nuovo Cimento {\bf 107B} (1992) pp. 703-711.
\bibitem{Cam95} Campos, R.G.; {\sl A Quadrature Formula for the Hankel Transform}, Numerical Algorithms, {\bf 9} (1995) pp. 343-354
\bibitem{Cam08} Campos, R.G.; Dom\'{\i}nguez-Mota F.; Coronado, E., {\sl Quadrature formulas for integrals transforms generated by orthogonal polynomials},  (to appear in IMA J. Numer. Anal.) http://arxiv.org/abs/0805.2111v1 [math.NA].
\bibitem{Bul08} Bultheel, A.; Sulbaran, H.M., {\sl http://www.cs.kuleuven.ac.be/~nalag/research/software/FRFT/}
\bibitem{Sze75} Szeg\"o G., {\sl Orthogonal Polynomials}, Colloquium Publications, American Mathematical Society, Providence, Rhode Island, 1975.
\bibitem{Gol69} Golub G. H.; Welsch, J. H., {\sl Calculation of Gauss quadrature rules},  Math. Comp. {\bf 23} (1969) 221-230.
\bibitem{Gau99} Gautschi, W., {\sl Orthogonal polynomials and quadrature}, Electron. Trans. Numer. Anal. {\bf 9} (1999) 65-76.
\bibitem{Erd53} Erd\'elyi, A. (ed.), {\sl Higher Transcendental Functions}, Vols I and II, McGraw Hill, New York, 1953.
\bibitem{Pel06} Pellat-Finet P.; Fogret, E., {\sl Complex order fractional Fourier transforms and their use in diffraction theory. Application to optical resonators}, Optics Comm. {\bf 258} (2006) 103-113.
\bibitem{Sol09} Sol\'{\i}s-Ort\'{\i}z, S.; Campos, R.G.; F\'elix J.; Obreg\'on, O., {\sl Coincident frequencies and relative phases among brain activity and hormonal signals}, Behav. Brain Funct.  2009, {\bf 5}:18  doi:10.1186/1744-9081-5-18
\bibitem{Erd53b} Erd\'elyi, A. (ed.), {\sl Tables of Integral Transforms}, Vol 1, McGraw Hill, New York, 1953.
\bibitem{Obe73} Oberhettinger, F., {\sl Fourier Expansions, a Collection of Formulas}, Academic Press, New York, USA, 1973

\end{thebibliography}
\end{document}